%% file: inter_deriv.tex
\documentclass[review]{elsarticle}
\usepackage{lineno,hyperref}
\modulolinenumbers[5]

\newtheorem{theo}{Theorem}
\newtheorem{ftheo}{Inaccurate Theorem}
\newtheorem{coro}[theo]{Corollary}

\newtheorem{lemma}[theo]{Lemma}

\input{commands}

\journal{Information Sciences}

\begin{document}

\begin{frontmatter}

\title{On the differentiability of interval functions}

\author{Walter F. Mascarenhas}
\address{Universidade de S\~{a}o Paulo, Brazil}
\ead[url]{www.ime.usp.br/~walterfm}
\ead{walterfm@ime.usp.br}

\address[mymainaddress]{Rua do Mat\~{a}o 1010, S\~{a}o Paulo, SP, Brazil, CEP 05508-090}

\begin{abstract}
Two articles published by Information Sciences discuss the
derivatives of interval functions in the sense of Svetoslav Markov. 
The authors of these articles tried to characterize for which functions 
and points such derivatives exist.
Unfortunately, their characterization is inaccurate. This article
describes this inaccuracy and explains how it can be corrected.
\end{abstract}

\begin{keyword}
\texttt{Interval functions} \sep Derivatives
\end{keyword}

\end{frontmatter}

\linenumbers

\section{Introduction}
\label{intro}

In his 1979 article \cite{Markov}, 
Markov presented a Calculus for {\it interval functions} of a real variable,
i.e., functions $F$ of the form
$\wfc{F}{t} = [\wfc{f}{t},\wfc{g}{t}]$ for
$f,g: \Omega \subset \wrone{} \rightarrow \wrone{}$ with $\wfc{f}{t} \leq \wfc{g}{t}$ for
$t \in \Omega$. ($\wrone{}$ is the set of real numbers.)
Assuming that $x$ is in the interior of $\Omega$, as we
will ways do, he uses the operator $\ominus$, defined by
\[
[\underline{a}, \overline{a}] \ominus [\underline{b}, \overline{b}]
= \left[\, \wfc{\min}{\underline{a} - \underline{b}, \overline{a} - \overline{b}},
	 \wfc{\max}{\underline{a} - \underline{b}, \overline{a} - \overline{b}}\right],
\]
the usual division operator $/$ for intervals and numbers, and the Hausdorff metric
\[
\wfc{\wrm{dist}}{ [\underline{a},\overline{a}], [\underline{b},\overline{b}] } =
\wfc{\max}{\,
\wabs{\, \underline{a} - \underline{b} \, }, \,  
\wabs{\,  \overline{a} - \overline{b}\,}  \,
},
\]
to define the derivative of $F$ at $x$ as 
\begin{equation}
\label{md}
\wdfc{F}{x} := \lim_{t \rightarrow x} \frac{\wfc{F}{t} \ominus \wfc{F}{x}}{t - x}.
\end{equation}
The relevance of Markov's $\wdfc{F}{x}$, in theory and in practice, was latter discussed in 
\cite{Osuna1,Osuna2,Stefanini}, and readers unaware of it should read these articles 
and some of their references. 

Markov showed that if $f$ and $g$ are differentiable at $x$ then 
\[
\wdfc{F}{x} = \left[  \,
\wfc{\min}{\wdfc{f}{x},\wdfc{g}{x}}, \,
\wfc{\max}{\wdfc{f}{x},\wdfc{g}{x}}  \, 
\right],
\]
and he also showed that $\wdfc{F}{x}$ may exist in cases in which 
$\wdfc{f}{x}$ or $\wdfc{g}{x}$ does not exist. To illustrate this fact,
in page 331 of \cite{Markov}, he considered the directional derivatives
\[
\wdfcs{f}{-}{x} = \lim_{t \uparrow 0} \frac{\wfc{f}{t} - \wfc{f}{x}}{t - x}
\hspace{1cm} \wrm{and} \hspace{1cm}
\wdfcs{f}{+}{x} = \lim_{t \downarrow 0} \frac{\wfc{f}{t} - \wfc{f}{x}}{t - x}
\]
and claimed that if $\wdfcs{f}{\pm}{x}$ and $\wdfcs{g}{\pm}{x}$ 
exist then $\wdfc{F}{x}$ exists when
\begin{equation}
\label{dpm}
\wdfcs{f}{-}{x} = \wdfcs{g}{+}{x}
\whs{1} \wrm{and} \whs{1}
\wdfcs{g}{-}{x} = \wdfcs{f}{+}{x}.
\end{equation}

The authors of \cite{Osuna1} and \cite{Osuna2} misunderstood page 331
of \cite{Markov} and this misunderstanding lead them to ``prove'' 
a ``theorem'' which in our notation reads like

\begin{ftheo} 
\label{wrong}
(inaccurate) Theorem 6 in \cite{Osuna1}, (inaccurate) Theorem 1 in \cite{Osuna2}.
The function $\wfc{F}{t} = [\wfc{f}{t}, \wfc{g}{t}]$  
has a derivative $\wdfc{F}{x}$ at $x$ if AND ONLY IF one of the following
cases hold:
\begin{itemize}
\item[(a)] The functions $f$ and $g$ are differentiable at $x$.
\item[(b)] The derivatives $\wdfcs{f}{\pm}{x}$ and 
$\wdfcs{g}{\pm}{x}$ 
exist and Equation \wref{dpm} holds.\wqes{}
\end{itemize}
\end{ftheo}

Markov did not state the ``ONLY IF'' part of this  ``theorem.'' In fact,
the next lemma shows that this part of ``Theorem 1'' is false:

\begin{lemma}
\label{counter}
Let $\wq{}$ be the set of rational numbers.
The function $\wfc{F}{t} = [\wfc{f}{t}, \wfc{g}{t}]$ for
\[
\wfc{f}{t} := 
\left\{
\begin{array}{r}
t \ \wrm{if}  \ t \in \wq{} \cap \wlr{-1,1}, \\ 
0 \ \wrm{if}  \ t \in \wlr{-1,1} \setminus{} \wq{}. 
\end{array}
\right.
\whs{1.0} \wrm{and} \whs{1.0}
\wfc{g}{t} := 
\left\{
\begin{array}{r}
1 \ \wrm{if}      \ t \in \wq{} \cap \wlr{-1,1},  \\
t + 1 \ \wrm{if}  \ t \in (-1,1) \setminus{} \wq{}.
\end{array}
\right.
\]
is such that $\wdfc{F}{0} = [0,1]$, but
$\wdfcs{f}{\pm}{0}$ and $\wdfs{g}{\pm}{0}$ do not exist.\wqes{}
\end{lemma}

The ``IF'' part of Theorem \ref{wrong}, which is correct, follows from 
the next theorem:

\begin{theo}
\label{thmif} 
If the right derivatives $\wdfcs{f}{+}{x}$ and $\wdfcs{g}{+}{x}$ exist
then the right derivative
\[
\wdfcs{F}{+}{x} := \lim_{t \downarrow 0} \frac{\wfc{F}{t} \ominus \wfc{F}{x}}{t - x}
\]
exists and 
\begin{equation}
\label{rightd}
\wdfcs{F}{+}{x} = [\wfc{\min}{ \wdfcs{f}{+}{x}, \wdfcs{g}{+}{x}}, 
                  \wfc{\max}{ \wdfcs{f}{+}{x}, \wdfcs{g}{+}{x}} ].
\end{equation}
Similarly, if the left derivatives 
$\wdfcs{f}{-}{x}$ and $\wdfcs{g}{-}{x}$ exist 
then the left derivative
\[
\wdfcs{F}{-}{x} := \lim_{t \uparrow 0} \frac{\wfc{F}{t} \ominus \wfc{F}{x}}{t - x}
\]
exists and 
\begin{equation}
\label{leftd}
\wdfcs{F}{-}{x} = [\wfc{\min}{ \wdfcs{f}{-}{x}, \wdfcs{g}{-}{x}}, 
                  \wfc{\max}{ \wdfcs{f}{-}{x}, \wdfcs{g}{-}{x}} ].\wqes{}
\end{equation}
\end{theo}

Theorem \ref{thmif} can be proved with the same 
technique used by Markov to prove his Theorem 3. 
We omit this proof for the sake of brevity, and because our
goal is to provide conditions under which we can fix the 
``ONLY IF'' part of ``Theorem 1.''

The most general condition on $f$ and $g$ which we could find in order
to fix ``Theorem 1'' is the existence of
functions $\alpha,\beta,c,d: \Omega \rightarrow \wrone{}$ 
such that 
$\alpha, \beta$ and $c$ are continuous in $\Omega$, $\wfc{\alpha}{t} \neq \wfc{\beta}{t}$ for $t \in \Omega$,
$\alpha$ and $\beta$ are bounded, $d$ is differentiable at $x$, with $\wfc{d}{x} = \wdfc{d}{x} = 0$ and
\begin{equation}
\label{cond}
\wfc{\alpha}{t} \wfc{f}{t} + \wfc{\beta}{t} \wfc{g}{t} = \wfc{c}{t} + \wfc{d}{t}.
\end{equation}
This condition is reasonable and applies to most practical situations. 
For instance, it holds if $f$ is continuous in $\Omega$, in which case we can take $\alpha = 1$, $\beta = 0$, $c = f$ and $d = 0$.
It also holds when the length of $\wfc{F}{t}$, $\wfc{g}{t} - \wfc{f}{t}$, is continuous, in which case we can take
$\alpha = -1$, $\beta = 1$, $c = g - f$ and $d = 0$. 
We formalize these results in the next theorem, which is
proved with other results in Section Proofs.

\begin{theo}
\label{thm} 
Consider functions $\alpha,\beta, c, d: \Omega \rightarrow \wrone{}$ 
and $\mu > 0$ such that $\wfc{d}{x} = 0$ and 
$ \wabs{\wfc{\alpha}{t}} + \wabs{\wfc{\beta}{t}} \leq \mu$ 
and $\mu \wabs{\, \wfc{\alpha}{t} - \wfc{\beta}{t} \,} \geq 1$ for all $t \in \Omega$.
If $\wdfcs{F}{+}{x}$
exists, Equation \wref{cond} holds for 
$t \in \Omega \cap (x,+\infty)$, 
$\alpha, \beta$ and $c$ are
 continuous in $\Omega \cap (x,+\infty)$ and $\wdfcs{d}{+}{0} = 0$,
then $\wdfcs{f}{+}{x}$ and $\wdfcs{g}{+}{x}$ exist.
If $\wdfcs{F}{-}{x}$
exists, Equation \wref{cond} holds for $t \in \Omega \cap (-\infty,x)$, 
$\alpha, \beta$ and $c$ are continuous in $\Omega \cap (-\infty,x)$ and $\wdfcs{d}{-}{0} = 0$,
then $\wdfcs{f}{-}{x}$ and $\wdfcs{g}{-}{x}$ exist.\wqes{}
\end{theo}

The last two theorems lead to  a proper ``ONLY IF'' part for ``Theorem 1'':
\begin{coro}
If $\wdfc{F}{x}$ exists and there exist functions $\alpha,\beta, c,d: \Omega \rightarrow \wrone{}$ 
such that $\alpha$, $\beta$ and $c$ are continuous, 
$\wfc{\alpha}{t} \neq \wfc{\beta}{t}$ for $t \in \Omega$, 
$\wfc{d}{x} = \wdfc{d}{x} = 0$ and
$\alpha f + \beta g = c + d$ 
then the directional derivatives $\wdfcs{f}{\pm}{x}$ and 
$\wdfcs{g}{\pm}{x}$  exist and
\begin{eqnarray}
\nonumber
\wdfc{F}{x} & = &
\left[\, \wfc{\min}{ \wdfcs{f}{-}{x}, \wdfcs{g}{-}{x}}, \,
 \wfc{\max}{ \wdfcs{f}{-}{x}, \wdfcs{g}{-}{x}} \, \right] \\
\label{ufa}
&  =  & \left[\, \wfc{\min}{ \wdfcs{f}{+}{x}, \wdfcs{g}{+}{x}}, \, 
 \wfc{\max}{ \wdfcs{f}{+}{x}, \wdfcs{g}{+}{x}} \, \right].
\end{eqnarray}
In particular, if
$\wdfc{F}{x}$ exists and $f$ or $g$ or $g - f$ is continuous, then 
the directional derivatives $\wdfcs{f}{\pm}{x}$ and $\wdfcs{g}{\pm}{x}$ 
exist and Equation \wref{ufa} holds. \wqes{}
\end{coro}

Combining this corollary with the next lemma, we end this introduction 
with the simplest version of the ONLY IF part of ``Theorem 1'' that we could find:
\begin{lemma}
\label{lemCont}
If $\wdfc{F}{x}$ exists then $\wfc{f}{t}$ and $\wfc{g}{t}$ are continuous at $x$.\wqes{}
\end{lemma}

\begin{coro}
\label{simple}
If $\wdfc{F}{x}$ exists for all $x$ in an open interval $(a,b)$ 
then the directional derivatives $\wdfcs{f}{\pm}{x}$ and 
$\wdfcs{g}{\pm}{x}$  exist for all $x \in (a,b)$ and satisfy Equation \wref{ufa}. 
\wqes{}
\end{coro}
The Denjoy-Young-Saks family of theorems
show that, in Corollary \ref{simple}, 
the set $\wcal{S}$ of $x \in (a,b)$ at which
$f$ or $g$ is not differentiable is ``tiny.'' 
Formally, the Lebesgue measure of $\wcal{S}$  is zero and it
is of first category
(in Baire's sense.) Actually, it is $\sigma$-porous (see \cite{Zajicek}.)

Thanks to the results above, it is quite likely that if 
the inaccurate
``ONLY IF'' part of ``Theorem 1'' was ever used then it caused no harm,
because their hypothesis would hold and ensure the validity
of the conclusions derived from ``Theorem 1.'' However, we still
believe that the minor technical point discussed in this article deserves
the attention of the potential users of ``Theorem 1.''

\section{Proofs}

{\bf Proof of Lemma \ref{counter}}
We show that 
\begin{equation}
\wfc{q}{t} =  \frac{ \wfc{F}{t} \ominus{} \wfc{F}{0} }{t} = [0,1]
\end{equation}
for $t \in (-1,1) \setminus \wset{0}$, by considering the four cases below:
\begin{itemize}
\item if $t > 0$ and $t \in \wq{}$ then
\[
\wfc{q}{t} =
\frac{
[\wfc{f}{t}, \wfc{g}{t}] \ominus [\wfc{f}{0}, \wfc{g}{0}]}{t}
= \frac{[t, 1]\ominus [0,1]}{t}  
\]
\[
= \frac{[ \wfc{\min}{t - 0, 1 - 1}, \wfc{\max}{t - 0, 1 - 1} ]}{t}  = \frac{[0,t]}{t}  = [0,1].
\]
\item if $t < 0$ and $t \in \wq{}$ then
\[
\wfc{q}{t} =
\frac{[\wfc{f}{t}, \wfc{g}{t}] \ominus [\wfc{f}{0}, \wfc{g}{0}]}{t}
= \frac{[t, 1] \ominus [0,1] }{t}  
\]
\[
= \frac{[ \wfc{\min}{t - 0, 1 - 1}, \wfc{\max}{t - 0, 1 - 1} ]}{t} 
 = \frac{[t,0]}{t}  = [0,1].
\]
\item if $t > 0$ and $t \not{\in} \wq{}$ then
\[
\frac{ [\wfc{f}{t}, \wfc{g}{t}] \ominus [\wfc{f}{0}, \wfc{g}{0}]}{t}
= \frac{[0, t + 1] \ominus [0,1] }{t} 
\]
\[
= \frac{[ \wfc{\min}{0 - 0, t + 1 - 1}, \wfc{\max}{0 - 0, t + 1 - 1} ]}{t} 
 =\frac{[0,t]}{t} = [0,1].
\]
\item if $t < 0$ and $t \not{\in} \wq{}$ then
\[
\wfc{q}{t} =
\frac{
[\wfc{f}{t}, \wfc{g}{t}] \ominus [\wfc{f}{0}, \wfc{g}{0}]
}{t}
= \frac{[0, t + 1] \ominus [0,1] }{t}  
\]
\[
= \frac{[ \wfc{\min}{0 - 0, t + 1 - 1}, 
          \wfc{\max}{0 - 0, t + 1 - 1} ]}{t} 
 = \frac{[t,0]}{t}  = [0,1]. \wqed{}
\]
\end{itemize}

{\bf Proof of Lemma \ref{lemCont}}
We show that if $f$ is discontinuous at $x$ then $\wdfc{F}{x}$ does not
exist (and the same holds for $g$.) 
In fact, in this case there exists a sequence $t_k$ which converges
to $x$ such that $\wfc{f}{t_k}$ does not converge to $\wfc{f}{x}$. This
implies that 
\[
\lim_{k \rightarrow \infty} \wabs{\frac{\wfc{f}{t_k} - \wfc{f}{x}}{t_k - x}}  = +\infty
\]
and the limit in the definition \wref{md} of $\wdfc{F}{x}$ does not exist. \wqed{}

{\bf Proof of Theorem \ref{thm}.}  
We only present the proof for the right derivatives. 
Since $x$ is in the interior of $\Omega$, there exists
$\delta > 0$ such that $[x, x + \delta) \subset \Omega$.
Writing $\wdfcs{F}{+}{x}$ as $[\underline{a},\overline{a}]$ we have that
\begin{equation}
\label{equd}
\underline{a} = \lim_{t \downarrow x} \frac{
\wfc{\min}{\wfc{f}{t} - \wfc{f}{x}, \, \wfc{g}{t} - \wfc{g}{x}}}{t - x}
\end{equation}
and
\begin{equation}
\label{eqod}
\overline{a} = \lim_{t \downarrow x} \frac{
\wfc{\max}{\wfc{f}{t} - \wfc{f}{x}, \, \wfc{g}{t} - \wfc{g}{x}}}{t - x}.
\end{equation}

Let us then define the set
\[
\wcal{P} := \wset{t \in (x,x + \delta) \ \wrm{with} \
\wfc{f}{t} - \wfc{f}{x} \leq \wfc{g}{t} - \wfc{g}{x}}.
\]
If $\wcal{P}$ is empty, let us define $i_p := x + \delta$. Otherwise, let $i_p$ be its infimum.
If $i_p > x$ then 
\[
\frac{\wfc{\min}{\wfc{f}{t} - \wfc{f}{x}, \, \wfc{g}{t} - \wfc{g}{x}}}{t - x} = 
\frac{\wfc{g}{t} - \wfc{g}{x}}{t - x}
\]
for $t \in (x,i_p)$ and Equation
\wref{equd} implies that $\wdfcs{g}{+}{x}$ exists and is equal to $\underline{a}$.
Similarly, for $t \in (x,i_p)$,
\[
\frac{\wfc{\max}{\wfc{f}{t} - \wfc{f}{x}, \, \wfc{g}{t} - \wfc{g}{x}}}{t - x} = 
\frac{\wfc{f}{t} - \wfc{f}{x}}{t - x}
\]
and Equation \wref{eqod} implies that 
$\wdfcs{f}{+}{x} = \overline{a}$. Therefore, if $i_p > x$ then
we are done.

Analogously, if the set
\[
\wcal{Q} := \wset{t \in (x,x + \delta) \ \wrm{with} \
\wfc{f}{t} - \wfc{f}{x} \geq \wfc{g}{t} - \wfc{g}{x}}.
\]
is empty or if its infimum $i_q$ is greater than $x$ then the 
left derivatives $\wdfcs{f}{+}{x}$ and $\wdfcs{g}{+}{x}$ exist,
and we are done.
Therefore, we can assume that $\wcal{P}$ and $\wcal{Q}$ are not empty and
$i_p = i_q = x$. We claim that in this case $\underline{a} = \overline{a}$,
and prove this claim by showing that
\begin{equation}
\label{tiny_da}
\overline{a} - \underline{a} = \wabs{ \, \overline{a} - \underline{a} \,} \leq 2 \mu \wlr{1 + \mu} \epsilon
\end{equation}
for all $\epsilon > 0$.
Let us then prove Equation \wref{tiny_da}. Equations \wref{equd} and \wref{eqod} 
and the hypothesis that $\wfc{d}{x} = \wdfcs{d}{+}{x} = 0$
imply that there exists $t_{\epsilon} \in (x, x+ \delta)$ such that if 
$t \in (x,t_{\epsilon})$ then
\begin{eqnarray}
\label{bnd}
\wabs{\wfc{d}{t}} \leq \epsilon \wlr{t - x}, \\
\label{bndp}
\wabs{\, \frac{
\wfc{\min}{\wfc{f}{t} - \wfc{f}{x}, \, \wfc{g}{t} - \wfc{g}{x}}}{t - x}
- \underline{a}  
\,}
\leq \epsilon, \\
\label{bndq}
\wabs{\, \frac{
\wfc{\max}{\wfc{f}{t} - \wfc{f}{x}, \, \wfc{g}{t} - \wfc{g}{x}}}{t - x}
- \overline{a}  
\,}
\leq \epsilon.
\end{eqnarray}
Consider 
$\wcal{P}_\epsilon := \wcal{P} \cap (x, t_\epsilon)$
and $\wcal{Q}_\epsilon := \wcal{Q} \cap (x, t_\epsilon)$.
These sets are not empty because $\inf \wcal{P} = \inf \wcal{Q} = x$.
Calling by $\overline{\wcal{P}}_\epsilon$ and  $\overline{\wcal{Q}}_\epsilon$
their closures in $(x,t_\epsilon)$, we have that 
$(x,t_\epsilon) = \wcal{P}_\epsilon \cup \wcal{Q}_\epsilon = 
\overline{\wcal{P}}_\epsilon \cup  \overline{\wcal{Q}}_\epsilon$.
Since $(x,t_\epsilon)$ is connected, 
$\overline{\wcal{P}}_\epsilon \cap  \overline{\wcal{Q}}_\epsilon$
is not empty. Taking $t \in \overline{\wcal{P}}_\epsilon \cap  \overline{\wcal{Q}}_\epsilon$,
there exist sequences $p_k \in \wcal{P}_\epsilon$ and $q_k \in \wcal{Q}_\epsilon$
such that $\lim_{k\rightarrow \infty} p_k = \lim_{k \rightarrow \infty} q_k = t$.
The definitions of $\wcal{P}$ and $\wcal{Q}$ and Equations \wref{bnd} and \wref{bndq} yield
\begin{eqnarray}
\nonumber
\wabs{\, \wfc{d}{p_k} \, } \leq \epsilon \wlr{p_k - x}
\whs{1} \wrm{and} \whs{1}
\wabs{\,  \wfc{d}{q_k}\, } \leq \epsilon \wlr{q_k - x},\\
\nonumber
\wabs{\, \frac{
\wfc{f}{p_k} - \wfc{f}{x}}{p_k - x} - \underline{a}  \,}
\leq \epsilon
\whs{0.5} \wrm{and} \whs{0.5}
\wabs{\, \frac{
\wfc{f}{q_k} - \wfc{f}{x}}{q_k - x}  - \overline{a}  \,}
\leq  \epsilon, \\
\nonumber
\wabs{\, \frac{
\wfc{g}{q_k} - \wfc{g}{x}}{q_k - x} - \underline{a}  \,}
\leq  \epsilon
\whs{0.5} \wrm{and} \whs{0.5}
\wabs{\, \frac{
\wfc{g}{p_k} - \wfc{g}{x}}{p_k - x}  - \overline{a}  \,}
\leq \epsilon.
\end{eqnarray}
These equations imply that the sequences $\wfc{f}{p_k}$, 
$\wfc{f}{q_k}$, $\wfc{g}{p_k}$, $\wfc{g}{q_k}$, $\wfc{d}{p_k}$ and $\wfc{d}{q_k}$ are bounded
and, by taking subsequences if necessary, we can assume that
they have limits $f_p, f_q, g_p$, $g_q$, $d_p$ and $d_q$ respectively.
We also have that
\begin{eqnarray}
\nonumber
\wfc{\alpha}{p_k} \wfc{f}{p_k} + \wfc{\beta}{p_k} \wfc{g}{p_k} = \wfc{c}{p_k} + \wfc{d}{p_k}, \\
\nonumber
\wfc{\alpha}{p_k} \wfc{f}{q_k} + \wfc{\beta}{p_k} \wfc{g}{q_k} = \wfc{c}{q_k} + \wfc{d}{q_k}.
\end{eqnarray}
Taking the limit $k \rightarrow \infty$ in the equations above, 
and recalling that $\alpha, \beta$ and $c$ are continuous, we get
\begin{eqnarray}
\label{dd}
\wabs{\, d_p  \, } \leq \epsilon \wlr{t - x}
\whs{1.05} \wrm{and} \whs{0.6}
\wabs{\,  d_q \, } \leq \epsilon \wlr{t - x}, \whs{0.6} \\
\label{d3}
\wabs{\, \frac{
f_p - \wfc{f}{x}}{t - x} - \underline{a}  \,}
\leq \epsilon
\whs{0.5} \wrm{and} \whs{0.5}
\wabs{\, \frac{
f_q - \wfc{f}{x}}{t - x}  - \overline{a}  \,}
\leq \epsilon, \\
\label{d4}
\wabs{\, \frac{
g_q - \wfc{g}{x}}{t - x} - \underline{a}  \,}
\leq \epsilon
\whs{0.5} \wrm{and} \whs{0.5}
\wabs{\, \frac{
g_p - \wfc{g}{x}}{t - x}  - \overline{a}  \,}
\leq \epsilon.\\
\label{dc}
\wfc{\alpha}{t} f_p + \wfc{\beta}{t} g_p = \wfc{c}{t} + d_p
\whs{0.5} \wrm{and} \whs{0.5}
\wfc{\alpha}{t} f_q + \wfc{\beta}{t} g_q = \wfc{c}{t} + d_q.
\end{eqnarray}
Combining equations \wref{dd} and \wref{dc} we obtain
that
\begin{equation}
\label{foo}
\wabs{ \,
\wfc{\alpha}{t} \frac{f_p - f_q}{t - x}  + 
\wfc{\beta}{t}  \frac{g_p - g_q}{t - x} \, } \leq 2 \epsilon.
\end{equation}
Equations \wref{d3} and \wref{d4} yield
\[
\wabs{\, \frac{f_p - f_q}{t - x} + \wlr{ \overline{a} - \underline{a}}  \,}
\leq 2 \epsilon
\whs{0.5} \wrm{and} \whs{0.5}
\wabs{\, \frac{
g_p - g_q}{t - x} - \wlr{\overline{a} - \underline{a}}  \,}
\leq 2 \epsilon,
\]
and these inequalities yield
\begin{equation}
\label{bar}
\wabs{\, \wfc{\alpha}{t} \frac{f_p - f_q}{t - x} +
         \wfc{\beta}{t} \frac{g_p - g_q}{t - x} 
+ \wlr{\wfc{\alpha}{t} - \wfc{\beta}{t}} \wlr{\overline{a} - \underline{a}}  \,}
\leq 2 \wlr{\wabs{\wfc{\alpha}{t}} + \wabs{\wfc{\beta}{t}} } \epsilon.
\end{equation}
Combining inequalities \wref{foo} and \wref{bar}
we complete the proof of Equation \wref{tiny_da}.

We have now shown that $\underline{a} = \overline{a}$,
and define $a := \underline{a} = \overline{a}$.
Let us now prove that $\wdfcs{f}{+}{x} = a$.
In fact, given $\epsilon > 0$, Equations \wref{equd} and \wref{eqod} imply that there exists 
$t_\epsilon \in (x, x + \delta)$ 
such that if $t \in (x,t_\epsilon)$ then
\[
-\epsilon \leq \frac{\wfc{\min}{\wfc{f}{t} - \wfc{f}{x}, 
                  \wfc{g}{t} - \wfc{g}{x}}}{t - x} - a 
\leq \frac{\wfc{f}{t} - \wfc{f}{x}}{t - x} - a
\]
\[
\leq \frac{\wfc{\max}{\wfc{f}{t} - \wfc{f}{x}, 
              \wfc{g}{t} - \wfc{g}{x}}}{t - x} - a 
              \leq \epsilon.
\]
Therefore, if $t \in (x, t_\epsilon)$ then
\[
\wabs{\, \frac{\wfc{f}{t} - \wfc{f}{x}}{t - x} - a \, } \leq \epsilon,
\]
and this shows that $\wdfcs{f}{+}{x} = a$. Analogously,
$\wdfcs{g}{+}{x} = a$ and we are done.\wqed{}

\bibliography{inter_deriv}

\end{document}

%% file: commands.tex
\usepackage{type1cm}   
\usepackage{dsfont}         
\usepackage{makeidx}         
\usepackage{graphicx}        
\usepackage{multicol}        
\usepackage[bottom]{footmisc}

\usepackage{newtxtext}       %
\usepackage{newtxmath}       
\usepackage{pict2e}

\usepackage{algorithm}
\usepackage{algpseudocode}
\usepackage{amsmath}
\usepackage{subcaption}
\usepackage{mathrsfs}

\newcommand{\wabs}[1]{{\left| {#1} \right|}}

\newcommand{\wcal}[1]{{\cal {#1}}}

\newcommand{\wdfs}[3]{{#1}'_{#2} \left(#3\right)}
\newcommand{\wdfcs}[3]{{#1}'_{#2}\!\left(#3\right)}

\newcommand{\wdfc}[2]{{#1}'\!\left(#2\right)}

\newcommand{\wfc}[2]{{#1}\!\left(#2\right)}

\newcommand{\whs}[1]{\hspace{#1cm}}

\newcommand{\wq}{{\mathds Q}}

\newcommand{\wqed}{\hfill \ensuremath{\blacksquare}}
\newcommand{\wqes}{\hfill \ensuremath{\blacktriangle}}

\newcommand{\wlr}[1]{\left({#1}\right)}

\newcommand{\wref}[1]{$\wlr{\ref{#1}}$}

\newcommand{\wrone}{\mathds R}

\newcommand{\wrm}[1]{\mathrm{#1}}

\newcommand{\wset}[1]{\left\{ {#1} \right\}}
\newcommand{\wtr}[1]{\mathrm{T}}